\documentclass[10pt,a4paper]{article}
\usepackage[top=70pt,bottom=90pt,left=64pt,right=64pt]{geometry}
\usepackage{amsmath} 
\usepackage{amsthm}
\usepackage{amssymb} 
\usepackage{amsfonts}
\usepackage{mathrsfs}
\usepackage{graphicx} 
\usepackage{inputenc}
\usepackage{upref} 
\usepackage{float}
\usepackage[all]{xy}
\usepackage{bm}
\usepackage{calc}
\usepackage{hyperref}

\newtheorem{theorem}{Theorem}[section]

\newtheorem{lemma}[theorem]{Lemma}
\newtheorem{definition}[theorem]{Definition}
\newtheorem{construction}[theorem]{Construction}
\newtheorem{proposition}[theorem]{Proposition}

\newcommand{\Z}{\mathbb{Z}}

\title{Twisted Bar Construction}
\author{Yutao Liu \\ University of Chicago
\\ Email: \href{mailto:yutao492@math.uchicago.edu}{yutao492@math.uchicago.edu}}

\begin{document}
\maketitle

\paragraph{Abstract:} This is an expository paper about a $\mathcal{C}_2$-equivariant analog of the bar construction. We will explain its construction and prove how it acts as a $\sigma$-fold delooping machine on twisted monoids, where $\sigma$ denotes the sign representation.

\section{Introduction}

One of the important questions in algebraic topology is to study iterated loop spaces and delooping machines. One classical way to construct the delooping machine for a topological monoid $A$ is the bar construction $BA$. When $A$ is homotopy commutative, the natural map from $A$ to the loop space $\Omega BA$ is a group completion in the sense that the induced map on $\pi_0$ is a group completion, and for any field of coefficients $k$, the induced map on homology $H_*(A;k)[\pi_0(A;k)^{-1}]\rightarrow H_*(\Omega BA;k)$ is an isomorphism.

A more general way to construct this delooping machine is to use the two-sided bar construction $B(\Sigma,C_1,A)$ for any $E_1$-algebra $A$ (where $C_1$ is an $E_1$-monad), which is equivalent to $BA$ when $A$ is a topological monoid (see [Thomason79]). When $A$ is an $E_n$-algebra, the analogous construction $B(\Sigma^n,C_n,A)$ acts as the $n$-fold delooping machine, which is equivalent to the iterated bar construction $B^nA$ when $A$ is a commutative monoid (which is naturally an $E_\infty$-algebra).

Both constructions are generalized to the equivariant case. [GM17] constructs the $V$-fold delooping machine for any representation $V$ by applying the two-sided bar construction on $E_V$-algebras. It shows that, for any $E_V$-algebra $A$, where $V$ contains at least one copy of the trivial representation, there exists a natural map from $A$ to $\Omega^V B(\Sigma^V,C_V,A)$ which is an equivariant group completion, in the sense that its restrictions on all fixed point spaces are group completions.

The classical bar construction $BA$, which is relatively easier for computations, is also generalized to $\mathcal{C}_2$-equivariant homotopy theory by defining a ``twisted" $\mathcal{C}_2$-action. We will call it ``twisted bar construction". It has the following property:

\begin{theorem}
There is a functor $B^\sigma$ from the category of twisted monoids to the category of based $\mathcal{C}_2$-spaces, such that, there is a natural $\mathcal{C}_2$-map $A\rightarrow\Omega^\sigma B^\sigma A$, which is a non-equivariant group completion if $A$ is non-equivariantly homotopy commutative, and its restriction on the fixed point spaces is a homotopy equivalence.
\end{theorem}

The idea for such a functor goes back at least to [BF85], where $B^\sigma$ is defined on certain monoids related to simplicial Hermitian rings and used to study the extended Hermitian algebraic K-theory. Theorem 1.1, which explains its key property, was first proved in [Dotto12] and [Stiennon13], with two different methods, and was generalized in [Moi13] to the simplicial case with more discussions of homology. The functor $B^\sigma$ is now being used in many computations concerning real spaces and spectra. One of the motivations of the author comes from [HS19], where $B^\sigma$ is used to analyze the $E_\sigma$-algebra structures on certain real spectra.

Moreover, the same idea of adding a ``twisted" $\mathcal{C}_2$-action also works for the cyclic bar construction, which provides one definition for THR, the real topological Hochschild homology (see [DMPR17] section 2). Section 4 in [Dotto12] also explains how the $\mathcal{C}_2$-equivariant delooping machine acts on THR.
\\

In this paper, we will give an explicit construction of $B^\sigma$ and prove Theorem 1.1. Our proof is similar to that in [Dotto12], but we fill in elementary details.

We introduce twisted monoids (which are sometimes called monoids with anti-involutions) in Section 2. They are algebras over a certain $E_\sigma$-operad of $\mathcal{C}_2$-spaces, and we compare that operad to other $E_\sigma$-operads in Proposition 2.3. More detailed discussions of $E_\sigma$-algebras are given in [Hil17].

In section 3, we define twisted simplicial spaces (which are sometimes called real simplicial spaces) and the twisted bar construction, and we use them to prove Theorem 1.1. We defer the most technical proof, that of Theorem 3.4, to section 5. It is motivated by the idea of Segal edgewise subdivisions [Segal73].

Finally, we show the compatibility of the twisted bar construction with the traditional bar construction in section 4, and we use this to construct a $V$-fold delooping machine for any $\mathcal{C}_2$-representation $V$.

\paragraph{Acknowledgements:} Emanuele Dotto kindly told me the historical background of this topic immediately after my original posting of this paper. He and Kristian Moi encouraged me not to withdraw the paper since it highlights a simple construction and its fundamental properties. I am very grateful to them both. I also thank my advisor, Peter May, for helpful suggestions on both the math and the English.

\section{Twisted Monoids}

From now on, we will only work in the category of based $G$-spaces in the case that $G=\mathcal{C}_2$.

In this section, we will introduce twisted monoids, which are the inputs of the twisted bar construction. They can be viewed as $\mathcal{C}_2$-spaces with a product compatible with the $\mathcal{C}_2$-action in a twisted way.

\begin{definition}
A $\mathcal{C}_2$-space $A$ is called a \textbf{twisted monoid} if it is a topological monoid in the non-equivariant sense, with the product satisfying $\tau(xy)=\tau(y)\tau(x)$ (Here $\tau$ is the generator of $\mathcal{C}_2$). A morphism between twisted monoids is a $\mathcal{C}_2$-map which commutes with the product.
\end{definition}

Let $\Sigma$ be the non-equivariant operad with components $\Sigma(n)=\Sigma_n$ and the usual structure maps. Then we can define a $\mathcal{C}_2$-action on it such that the action of $\tau$ on $\Sigma_n$ is the left multiplication by the element sending $1,2,...,n$ to $n,n-1,...,1$.

We can check that these actions commute with the operadic structure maps. For any $u_i\in\Sigma_{s_i}$, $i=1,2,...,t$; and $v\in\Sigma_t$, the structure map
$$\Sigma_t\times\Sigma_{s_1}\times...\times\Sigma_{s_t}\rightarrow\Sigma_s$$
(where $s=s_1+s_2+...+s_t$) sends $v\times u_1\times...\times u_t$ to the element $u\in\Sigma_s$, such that, for any $1\leq i\leq t$ and
$$s_1+...+s_{i-1}<j\leq s_1+...+s_i,$$
we have 
$$u(j)=u_i(j-\sum_{k=1}^{i-1}s_k)+\sum_{k=1}^{v(i)-1}s_{v^{-1}(k)}.$$

Thus $\tau v\times\tau u_1\times...\times\tau u_t$ is sent to the element $u^\prime\in\Sigma_s$ with
$$u^\prime(j)=\tau u_i(j-\sum_{k=1}^{i-1}s_k)+\sum_{k=1}^{\tau v(i)-1}s_{(\tau v)^{-1}(k)}.$$

Notice that, for any $w\in\Sigma_l$, we have
$$\tau w(a)=l+1-w(a)$$
and
$$(\tau w)^{-1}(a)=w^{-1}\tau(a)=w^{-1}(l+1-a)$$
for $a=1,2,...,l$. Thus

$$u^\prime(j)=s_i+1-u_i(j-\sum_{k=1}^{i-1}s_k)+\sum_{k=1}^{t+2-v(i)}s_{v^{-1}(t+1-k)}.$$

So we have

$$u(j)+u^\prime(j)=s_i+1+\sum_{k=1}^{v(i)-1}s_{v^{-1}(k)}+\sum_{k=1}^{t+2-v(i)}s_{v^{-1}(t+1-k)}.$$

The last two sums cover all of $s_1,s_2,...,s_t$ except $s_i$. So we have $u(j)+u^\prime(j)=s+1$, and hence $u^\prime=\tau u$.
\\

Therefore, the actions by $\mathcal{C}_2$ make $\Sigma$ into a $\mathcal{C}_2$-operad. We call it $\hat{\Sigma}$.

\begin{lemma}
A $\hat{\Sigma}$-algebra is equivalent to a twisted monoid.
\end{lemma}

\paragraph{Proof:} Let $A$ be a $\hat{\Sigma}$-algebra. Then we have a $\mathcal{C}_2$-structure map
$$\hat{\Sigma}(2)\times_{\Sigma_2}(A\times A)\rightarrow A.$$
Notice that $\hat{\Sigma}(2)$ is $\Sigma_2$ together with a left $\mathcal{C}_2$-action such that $\tau$ exchanges its two elements. Thus the underlying space of the left hand side is just $A\times A$, but the $\mathcal{C}_2$-action becomes $\tau(x,y)=(\tau y,\tau x)$, for any $x,y\in A$. So we get a product on $A$ which satisfies $\tau(xy)=\tau(y)\tau(x)$. Moreover, this product is associative in the non-equivariant sense, since by forgetting the $\mathcal{C}_2$-actions, it is exactly the product on a $\Sigma$-algebra. Therefore, $A$ is a twisted monoid.

On the other hand, if we have a twisted monoid $A$ first, we can make it into a $\hat{\Sigma}$-algebra by defining structure maps
$$\hat{\Sigma}(n)\times A^n\rightarrow A$$
$$(u,a_1,a_2,...,a_n)\mapsto a_{u^{-1}(1)}a_{u^{-1}(2)}...a_{u^{-1}(n)}.$$
This is a $\mathcal{C}_2$-map, since the action of $\tau$ sends $(u,a_1,...,a_n)$ to $(\tau u,\tau a_1,...,\tau a_n)$, whose image becomes 
$$(\tau a_{(\tau u)^{-1}(1)})...(\tau a_{(\tau u)^{-1}(n)})=(\tau a_{u^{-1}(n)})...(\tau a_{u^{-1}(1)})$$
which is exactly $\tau(a_{u^{-1}(1)}...a_{u^{-1}(n)})$.

We can check that such maps are compatible with the operadic maps by forgetting the $\mathcal{C}_2$-actions (then it becomes the natural $\Sigma$-algebra structure on any topological monoid). $\Box$
\\

As operads, $\Sigma$ and $\hat{\Sigma}$ can be viewed as simplified versions of the non-equivariant little $1$-disk and the equivariant little $\sigma$-disk operads.

\begin{proposition}
Let $\mathscr{D}_1$ and $\mathscr{D}_\sigma$ be the non-equivariant little $1$-disk operad and the equivariant little $\sigma$-disk operad. Then there exists a non-equivariant weak equivalence $\mathscr{D}_1\rightarrow\Sigma$, which induces an equivariant weak equivalence $\mathscr{D}_\sigma\rightarrow\hat{\Sigma}$. Thus each twisted monoid is automatically an $E_\sigma$-algebra, and [BM03] shows that the category of twisted monoids is Quillen equivalent to the category of $E_\sigma$-algebras.
\end{proposition}

\paragraph{Proof:} Notice that, each $\mathscr{D}_1(n)$ is the space of linear embeddings from $n$ copies of $(-1,1)$ to one copy of $(-1,1)$, which doesn't necessarily preserve the order. So an element in it can be expressed as a sequence of $n$ disjoint open subintervals $(a_1,b_1),...,(a_n,b_n)$ of $(-1,1)$. Reorder these intervals by choosing the unique element $u\in\Sigma_n$ such that
$$a_{u(1)}<b_{u(1)}\leq a_{u(2)}<b_{u(2)}\leq...\leq a_{u(n)}<b_{u(n)}.$$
Then we can define a map from $\mathscr{D}_1(n)$ to $\Sigma_n$ sending it to $u^{-1}$. These maps extend to an operadic map from $\mathscr{D}_1$ to $\Sigma$.

In the equivariant case, the little $\sigma$-disk operad $\mathscr{D}_\sigma$ can be viewed as $\mathscr{D}_1$ with an additional $\mathcal{C}_2$-action, such that $\tau$ sends $(a_1,b_1),...,(a_n,b_n)$ to its central symmetrical embedding $(-b_n,-a_n),...,(-b_1,-a_1)$. Thus if the image of $(a_1,b_1),...,(a_n,b_n)$ in $\Sigma_n$ is $u$, the image of $\tau((a_1,b_1),...,(a_n,b_n))$ will become $u$ left multiplied by the element sending $1,2,...,n$ to $n,n-1,...,1$. Notice that this is exactly how we define the $\mathcal{C}_2$-action on $\Sigma$ and obtain $\hat{\Sigma}$. So we get a $\mathcal{C}_2$-operadic map from $\mathscr{D}_\sigma$ to $\hat{\Sigma}$. Thus each twisted monoid is automatically an $E_\sigma$-algebra. Moreover, each $\mathscr{D}_1(n)$ is the disjoint union of contractible components. The map 
$$\mathscr{D}_1(n)\rightarrow\Sigma_n$$
is a 1-1 correspondence between the components of $\mathscr{D}_1(n)$ and discrete points in $\Sigma_n$, thus is a homotopy equivalence.

Now we add $\mathcal{C}_2$-actions on both sides. When $n>1$, $\tau$ sends each component of $\mathscr{D}_1(n)$ to a different one (and each point in $\Sigma_n$ to a different point). So the map above becomes a $\mathcal{C}_2$-equivalence. When $n=1$, we can write an explicit $\mathcal{C}_2$-homotopy from $\mathscr{D}_\sigma(1)$ to a single point as
$$\mathscr{D}_\sigma(1)\times [0,1]\rightarrow\mathscr{D}_\sigma(1)$$
$$(a,b)\times t\mapsto(-t+(1-t)a,t+(1-t)b).$$
So the map
$$\mathscr{D}_\sigma(1)\rightarrow\hat{\Sigma}_1\cong *$$
is still an equivalence. Thus the operadic map from $\mathscr{D}_\sigma$ to $\hat{\Sigma}$ is a weak equivalence. $\Box$

\section{Twisted Bar Construction}

Now we can explain the twisted bar construction.

\begin{construction}
For any twisted monoid $A$, construct $B^\sigma_* A$ in the same way as the non-equivariant bar construction, such that $B^\sigma_nA=A^n$. However, we define the $\mathcal{C}_2$-action on each $A^n$ by 
$$\tau(a_1,a_2,...,a_n)=(\tau a_n,\tau a_{n-1},...,\tau a_1).$$ 
Then the $\mathcal{C}_2$-actions commute with the face and degeneracy maps as $\tau\circ s_i=s_{n-i}\circ\tau$ and $\tau\circ d_i=d_{n-i}\circ\tau$.

In addition, define a $\mathcal{C}_2$-action on each
$$\Delta^n=\{(v_0,v_1,...,v_n)\in[0,1]^{n+1}:\sum v_i=1\}$$
by $\tau(v_0,v_1,...,v_n)=(v_n,v_{n-1},...,v_0)$. Now when we take the geometric realization
$$\coprod_n A^n\times\Delta^n/\sim$$
all the relations commute with the $\mathcal{C}_2$-actions. This makes the geometric realization into a $\mathcal{C}_2$-space. We call it the \textbf{twisted bar construction} of $A$, denoted by $B^\sigma A$.
\end{construction}

We can also make this idea more general:

\begin{definition}
A \textbf{twisted simplicial space} is a simplicial space $S_*$, together with a $\mathcal{C}_2$-action on each $S_n$, such that $\tau\circ s_i=s_{n-i}\circ\tau$ and $\tau\circ d_i=d_{n-i}\circ\tau$. Its geometric realization $|S_*|$ is a $\mathcal{C}_2$-space constructed as above.
\end{definition}

Assume that $A$ is a twisted monoid which is homotopy commutative in the non-equivariant sense. Our first goal is to show that $B^\sigma A$ is the degree $\sigma$ delooping of $A$ as follows:

\begin{theorem} \label{main}
Let $\Omega^\sigma B^\sigma A=Map_*(S^\sigma,B^\sigma A)$ be the space of based non-equivariant maps with $\mathcal{C}_2$ acting by conjugation. Then there exists a natural $\mathcal{C}_2$-map $A\rightarrow\Omega^\sigma B^\sigma A$, which is a non-equivariant group completion if $A$ is non-equivariantly homotopy commutative, and its restriction on the fixed point spaces is a homotopy equivalence.
\end{theorem}

Notice that, without the $\mathcal{C}_2$-action, $\Omega^\sigma B^\sigma A$ is the same as $\Omega BA$. And we already have a natural map $A\rightarrow\Omega BA$ defined by
$$a\mapsto (t\mapsto |a\times (t,1-t)|)$$
where $|a\times (t,1-t)|$ is the image of $a\times (t,1-t)\in A\times\Delta^1$ in the geometric realization. It's not hard to check that, if we use $\Omega^\sigma B^\sigma A$ instead, this will become a $\mathcal{C}_2$-map. Moreover, since this map is a non-equivariant group completion, it suffices to show that, the restricted map between the fixed point spaces is an equivalence.

We have 
$$(\Omega^\sigma B^\sigma A)^{\mathcal{C}_2}=Map^{\mathcal{C}_2}_*(S^\sigma,B^\sigma A)=\{f:[0,1]\rightarrow B^\sigma A:f(0)=f(1)=x_0,f(1-t)=\tau f(t),\forall t\in[0,1]\}$$
$$=\{f:[0,1]\rightarrow B^\sigma A:f(0)=x_0,f(1)\in(B^\sigma A)^{\mathcal{C}_2}\}$$
where $x_0$ is the base point of $B^\sigma A$. The last equality holds because each $f$ is determined by its value on $[\frac{1}{2},1]$. 
Let 
$$\iota:(B^\sigma A)^{\mathcal{C}_2}\hookrightarrow B^\sigma A$$
by the inclusion. Then we can rewrite it as
$$(\Omega^\sigma B^\sigma A)^{\mathcal{C}_2}=\{(f,x):x\in (B^\sigma A)^{\mathcal{C}_2}, f:[0,1]\rightarrow B^\sigma A,f(0)=x_0,f(1)=\iota(x)\}$$

Notice that this is exactly the homotopy fiber $F\iota$. The restricted map $A^{\mathcal{C}_2}\rightarrow (\Omega^\sigma B^\sigma A)^{\mathcal{C}_2}=F\iota$ now sends each $a\in A^{\mathcal{C}_2}$ to $(f,x)$, such that $f(t)=|a\times (1-t/2,t/2)|$, $x=|a\times(1/2,1/2)|$. Use $\varphi:A^{\mathcal{C}_2}\rightarrow F\iota$ to denote this restricted map.

In fact, we have a more explicit description for the fixed point space of $B^\sigma A$:

\begin{theorem} \label{mainlemma}
$(B^\sigma A)^{\mathcal{C}_2}$ is a quotient of the disjoint union of $(B^\sigma_{2n+1}A\times\Delta^{2n+1})^{\mathcal{C}_2}$, which consists of points with the form
$$(a_1,...,a_n,a,\tau a_n,...,\tau a_1)\times(v_0,...,v_n,v_n,...,v_0)$$
where $a\in A^{\mathcal{C}_2}$, for $n=0,1,...$. Moreover, $(B^\sigma A)^{\mathcal{C}_2}$ is homeomorphic to the two-sided bar construction $B(*,A,A^{\mathcal{C}_2})$, where the action of $a\in A$ on $A^{\mathcal{C}_2}$ is given by $b\mapsto ab\tau(a)$. The homeomorphism 
$$\eta:(B^\sigma A)^{\mathcal{C}_2}\rightarrow B(*,A,A^{\mathcal{C}_2})$$
can be expressed as
$$(a_1,...,a_n,a,\tau a_n,...,\tau a_1)\times (v_0,...,v_n,v_n,...,v_0)\mapsto (a_1,...,a_n,a)\times (2v_0,...,2v_n)\in B_n(*,A,A^{\mathcal{C}_2})\times\Delta^n.$$

The inclusion 
$$\iota: (B^\sigma A)^{\mathcal{C}_2}\hookrightarrow B^\sigma A$$
is homotopic to the natural projection $$p:B(*,A,A^{\mathcal{C}_2})\rightarrow B(*,A,*)=BA$$
(notice that $BA$ is the same as $B^\sigma A$ non-equivariantly). We can define an explicit homotopy
$$H:(B^\sigma A)^{\mathcal{C}_2}\times I\rightarrow B^\sigma A$$
by
$$(a_1,...,a_n,a,\tau a_n,...,\tau a_1)\times (v_0,...,v_n,v_n,...,v_0)\times t\mapsto$$
$$(a_1,...,a_n,a,\tau a_n,...,\tau a_1)\times((1+t)v_0,...,(1+t)v_n,(1-t)v_n,...,(1-t)v_0).$$
Then we have $H(x,0)=\iota(x)$, $H(x,1)=p(\eta(x))$.
\end{theorem}

We will prove this theorem in section 5. For now, we will use this explicit homotopy to prove theorem 3.3.
\\

\paragraph{Proof of Theorem 3.3:}

Write $X=(B^\sigma A)^{\mathcal{C}_2}$, $Y=B(*,A,A^{\mathcal{C}_2})$. The homotopy $H:X\times I\rightarrow B^\sigma A$ in Theorem 3.4 induces a homotopy equivalence between the homotopy fibers of $\iota$ and $p$, as follows:

Notice that
$$F\iota=\{(f,x):x\in X, f:[0,1]\rightarrow B^\sigma A:f(0)=x_0,f(1)=\iota
(x)\},$$
$$Fp=\{(g,y):y\in Y, g:[0,1]\rightarrow B^\sigma A:g(0)=x_0,g(1)= p(y)\}.$$

So we can define a map $\alpha:F\iota\rightarrow Fp$ sending each $(f,x)\in F\iota$ to $(f\circ H(x,\cdot),\eta(x))$. Similarly, we have a map $\beta:Fp\rightarrow F\iota$ sending each $(g,y)\in Fp$ to $(g\circ H(\eta^{-1}(y),1-\cdot),\eta^{-1}(y))$. Both $\alpha\circ\beta$ and $\beta\circ\alpha$ are homotopic to the identity maps, since they just attach a loop, which is homotopic to the trivial one, at the end of each path. Thus they are homotopy inverses between $F\iota$ and $Fp$.

Consider the composed map
$$A^{\mathcal{C}_2}\xrightarrow{\varphi}F\iota\xrightarrow{\alpha}Fp.$$
It sends $a\in A^{\mathcal{C}_2}$ to $(f,x)$ with $f(t)=|a\times(1-t/2,t/2)|$, and $x=|a\times(1/2,1/2)|\in X$. Notice that $\eta(x)$ is the point $a\in A^{\mathcal{C}_2}=B_0(*,A,A^{\mathcal{C}_2})$. Thus $\alpha$ sends it further to $(g,a)$, where $g=f\circ H(|a\times(1/2,1/2)|,\cdot)$. So we have $g(t)=|a\times(1-t,t)|$ when $t\leq 1/2$, and $g(t)=|a\times(t,1-t)|$ when $t\geq 1/2$. Notice that such $g$ is homotopic to the trivial loop. Thus the composite map
$$\alpha\circ\varphi: A^{\mathcal{C}_2}\rightarrow Fp$$
is homotopic to the map sending $a\in A^{\mathcal{C}_2}$ to the pair of the trivial loop and the point $a$, which is exactly the inclusion from the actual fiber of $p$ (which is $A^{\mathcal{C}_2}$) to the homotopy fiber $Fp$. Therefore, both $\alpha\circ\varphi$ and $\alpha$ are homotopy equivalences. Thus $\varphi:A^{\mathcal{C}_2}\rightarrow (\Omega^\sigma B^\sigma A)^{\mathcal{C}_2}$ defined at the beginning is an equivalence. So the natural map $A\rightarrow \Omega^\sigma B^\sigma A$ is a non-equivariant group completion, and an equivalence when restricted on the fixed point spaces. $\Box$

\section{Deloopings of Commutative Twisted Monoids}

Now we can construct the $V$-fold delooping for any arbitrary representation $V$ with both the ordinary and the twisted bar constructions.

Notice that any representation $V$ can be written as $a+b\sigma$ with some $a,b\geq 0$. We can define the $V$-fold delooping as the composition of $a$ copies of $B$ and $b$ copies of $B^\sigma$. But we also have to show that this definition doesn't depend on the order of the composition.

\begin{proposition}
Let $S_*$, $T_*$ be twisted simplicial spaces. Then their product $S_*\times T_*$ (whose $n$-th component is $S_n\times T_n$) has a natural twisted simplicial space structure. Moreover, there exists a $\mathcal{C}_2$-homeomorphism between $|S_*\times T_*|$ and $|S_*|\times|T_*|$.
\end{proposition}

\paragraph{Proof:} The non-equivariant homeomorphism 
$$|S_*\times T_*|\rightarrow|S_*|\times|T_*|$$
is defined by sending the image of each
$$s_n\times t_n\times (v_0,v_1,...,v_n)\in S_n\times T_n\times\Delta^n$$
in the geometric realization to
$$|s_n\times (v_0,v_1,...,v_n)|\times|t_n\times (v_0,v_1,...,v_n)|\in|S_*|\times|T_*|.$$

The action by $\tau$ sends $s_n\times t_n\times (v_0,v_1,...,v_n)$ to $\tau s_n\times\tau t_n\times(v_n,v_{n-1},...,v_0)$, whose image becomes
$$|\tau s_n\times(v_n,v_{n-1},...,v_0)|\times|\tau t_n\times(v_n,v_{n-1},...,v_0)|\in|S_*|\times|T_*|,$$
which is exactly the image of $|s_n\times (v_0,v_1,...,v_n)|\times|t_n\times (v_0,v_1,...,v_n)|$ under the action of $\tau$. Thus the $\mathcal{C}_2$-action commutes with the homeomorphism, and makes it into an equivariant one. $\Box$
\\

Now we can talk about multi-simplicial objects which mix both the ordinary and the twisted objects.

Let $S_{*,*}$ be a bisimplicial space with a $\mathcal{C}_2$-action on each component, such that each $S_{p,*}$ becomes a simplicial object in the category of based $\mathcal{C}_2$-spaces, and each $S_{*,q}$ becomes a twisted simplicial space. Then we can define its geometric realization by realizing the first subscript (in the twisted way) and then the second one (in the ordinary way). We can also realize the second subscript first, which doesn't change the result, since both of them can be viewed as the quotient of
$$\coprod_{p,q}S_{p,q}\times\Delta^p\times\Delta^q$$
(where the action of $\mathcal{C}_2$ on $\Delta^p$ is trivial, but behaves as $\tau(v_0,...,v_q)=(v_q,...,v_0)$ on $\Delta^q$) by the relations on the first subscript and the relations on the second one. The assumptions on each $S_{p,*}$ and each $S_{*,q}$ guarantee that both kinds of relations are compatible with the $\mathcal{C}_2$-action.

In particular, we have the following property:

\begin{proposition}
Let $A$ be a commutative monoid in the category of based $\mathcal{C}_2$-spaces. Then both $BB^\sigma A$ and $B^\sigma BA$ are well-defined and have natural commutative monoid structures. Moreover, they are $\mathcal{C}_2$-homeomorphic to each other.
\end{proposition}

\paragraph{Proof:} Since a commutative monoid is both a monoid and a twisted monoid in the category of based $\mathcal{C}_2$-spaces, according to the commutativity between the product and the geometric realizations, both $B^\sigma A$ and $BA$ are still commutative monoids. Thus both $BB^\sigma A$ and $B^\sigma BA$ are well-defined as commutative monoids. Moreover, we can check that they are the two realizations (starting with different subscripts) of the same bisimplicial object as described above. So from the previous argument, they are $\mathcal{C}_2$-homeomorphic to each other. $\Box$
\\

Finally, we have the $V$-fold delooping as follows.

\begin{proposition}
For any commutative monoid $A$ in the category of based $\mathcal{C}_2$-spaces, the $V$-degree bar construction $B^VA$ is defined by applying the ordinary bar construction $a$ times and the twisted bar construction $b$ times for $V=a+b\sigma$. There exists a natural map $A\rightarrow\Omega^V B^VA$. When $A$ is $\mathcal{C}_2$-connected, this map is a $\mathcal{C}_2$-equivalence.
\end{proposition}

\section{Proof of Theorem 3.4}

We will divide the proof into several steps.

\paragraph{(i)} As the fixed point space, $(B^\sigma A)^{\mathcal{C}_2}$ is the quotient of $\coprod_n (A^n\times\Delta^n)^{\mathcal{C}_2}$ by the restricted relations on these components. A point $(a_1,a_2,...,a_n)\times(v_0,v_1,...,v_n)\in A^n\times\Delta^n$ is fixed if and only if
$$a_i=\tau a_{n+1-i},v_i=v_{n-i},i=n,n-1,...,\frac{n}{2}+1$$
when $n$ is even, or
$$a_i=\tau a_{n+1-i},i=n,n-1,...,\frac{n+3}{2}; a_{\frac{n+1}{2}}\in A^{\mathcal{C}_2}; v_i=v_{n-i},i=n,n-1,...,\frac{n+1}{2}$$
when $n$ is odd.

Notice that if $n$ is even and 
$$(a_1,a_2,...,a_n)\times(v_0,...,v_n)\in (A^n\times\Delta^n)^{\mathcal{C}_2},$$
then the degeneracy map $s_\frac{n}{2}$ induces an equivalence from it to 
$$(a_1,a_2,...,a_{\frac{n}{2}},e,a_{\frac{n}{2}+1},...,a_n)\times(v_0,v_1,...,v_{\frac{n}{2}-1},\frac{1}{2}v_{\frac{n}{2}},\frac{1}{2}v_{\frac{n}{2}},v_{\frac{n}{2}+1},...v_n)\in (A^{n+1}\times\Delta^{n+1})^{\mathcal{C}_2},$$
where $e$ is the unit of $A$ (which must be fixed by $\tau$). Thus $(B^\sigma A)^{\mathcal{C}_2}$ can be viewed as the quotient of the disjoint union of $(B^\sigma_{2n+1}A\times \Delta_{2n+1})^{\mathcal{C}_2}$ for $n=0,1,...$.

\paragraph{(ii)} Write $S_n=(A^{2n+1}\times\Delta^{2n+1})^{\mathcal{C}_2}$. Notice that each point in it is determined by $a_1,a_2,...,a_{n}\in A$, $a_{n+1}\in A^{\mathcal{C}_2}$, and $v_0,v_1,...,v_{n}$ whose sum is $1/2$. 

On the other hand, write 
$$T_n=B_n(*,A,A^{\mathcal{C}_2})\times\Delta^n=A^n\times A^{\mathcal{C}_2}\times\Delta^n.$$
Then we have a homeomorphism $T_n\rightarrow S_n$ defined by
$$(a_1,a_2,...,a_n)\times a\times(v_0,v_1,...,v_n)\mapsto(a_1,a_2,...,a_n,a,\tau a_n,...,\tau a_1)\times(\frac{1}{2}v_0,...\frac{1}{2}v_n,\frac{1}{2}v_n,...,\frac{1}{2}v_0)$$
(which is the inverse of the restriction of $\eta$ defined in Theorem \ref{mainlemma}).

Notice that $(B^\sigma A)^{\mathcal{C}_2}$ and $B(*,A,A^{\mathcal{C}_2})$ are quotients of $\coprod_n S_n$ and $\coprod_n T_n$ by some relations. It suffices to show that, if one relation from $S_n$ and $S_m$ identifies a pair of points $x\in S_n$ and $y\in S_m$, then their images in $T_n$ and $T_m$ are also identified by one or a sequence of relations, and vice versa.

\paragraph{(iii)} Any relation from $T_n$ to $T_m$ can be expressed by an order-preserving map 
$$s:\{0,1,...,m\}\rightarrow\{0,1,...,n\},$$
and it identifies 
$$(a_1,...,a_n)\times a\times (v_0,v_1,...,v_n)\in T_n$$ 
and 
$$(b_1,...,b_m)\times b\times (w_0,...,w_m)\in T_m,$$
if 
$$b_i=a_{s(i-1)+1}a_{s(i-1)+2}...a_{s(i)}, b=a_{s(m)+1}...a_n a \tau(a_n)...\tau(a_{s(m)+1}), v_j=\sum_{k\in s^{-1}(j)}w_k$$
for all $i,j$ (let $s(-1)=-1$). 

On the other hand, any relation from $S_n$ to $S_m$ can be expressed by an order-preserving map
$$t:\{0,1,...,2m+1\}\rightarrow\{0,1,...,2n+1\},$$
and it identifies
$$(c_1,c_2,...,c_{2n+1})\times(p_0,...,p_{2n+1})\in S_n\subset A^{2n+1}\times\Delta^{2n+1}$$
and 
$$(d_1,...,d_{2m+1})\times(q_0,...,q_{2m+1})\in S_m\subset A^{2m+1}\times\Delta^{2m+1},$$
(notice that some additional conditions are required on these $c_i,d_i,p_i,q_i$ since they are inside the fixed point space) if
$$d_i=c_{t(i-1)+1}c_{t(i-1)+2}...c_{t(i)}, p_j=\sum_{k\in t^{-1}(j)}q_k$$
for all $i,j$.

\paragraph{(iv)} One direction is easy to check. If we have an order preserving map
$$s:\{0,1,...,m\}\rightarrow\{0,1,...,n\},$$
identifying 
$$(a_1,...,a_n)\times a\times (v_0,v_1,...,v_n)\in T_n$$ 
and 
$$(b_1,...,b_m)\times b\times (w_0,...,w_m)\in T_m,$$
then their images in $S_n$ and $S_m$ are
$$(a_1,...,a_n,a,\tau a_n,...,\tau a_1)\times(\frac{1}{2}v_0,...\frac{1}{2}v_n,\frac{1}{2}v_n,...,\frac{1}{2}v_0)\in S_n$$
and
$$(b_1,...,b_m,b,\tau b_m,...,\tau b_1)\times(\frac{1}{2}w_0,...\frac{1}{2}w_m,\frac{1}{2}w_m,...,\frac{1}{2}w_0)\in S_m,$$
which are identified by the relation expressed by
$$t:\{0,1,...,2m+1\}\rightarrow\{0,1,...,2n+1\}$$
such that $t(i)=s(i)$ if $i\leq m$, and $t(i)=2n+1-s(2m+1-i)$ if $i\geq m+1$.

\paragraph{(v)} On the other hand, if a pair of points
$$(c_1,...,c_n,c,\tau c_n,...,\tau c_1)\times(p_0,...,p_n,p_n,...,p_0)\in S_n$$
and
$$(d_1,...,d_m,d,\tau d_m,...,\tau d_1)\times(q_0,...,q_m,q_m,...,q_0)\in S_m$$
is identified by an order-preserving map
$$t:\{0,1,...,2m+1\}\rightarrow\{0,1,...,2n+1\}$$
which is symmetric, in the sense that $t(i)+t(2m+1-i)=2n+1$ for all $i$, then their images
$$(c_1,...,c_n)\times c\times(2p_0,...,2p_n)\in T_n$$
and
$$(d_1,...,d_m)\times d\times(2q_0,...,2q_m)\in T_m$$
are identified by
$$s:\{0,1,..,m\}\rightarrow\{0,1,...,n\}$$
defined as the restriction of $t$. So it suffices to show that, any pair of equivalent points in $S_n$ and $S_m$ can be identified by one or a sequence of relations induced by symmetric maps.

\paragraph{(vi)} Assume that
$$x=(c_1,c_2,...,c_{2n+1})\times(p_0,...,p_{2n+1})=(c_1,...,c_n,c,\tau c_n,...,\tau c_1)\times(p_0,...,p_n,p_n,...,p_0)\in S_n$$
and
$$y=(d_1,...,d_{2m+1})\times(q_0,...,q_{2m+1})=(d_1,...,d_m,d,\tau d_m,...,\tau d_1)\times(q_0,...,q_m,q_m,...,q_0)\in S_m$$
are identified by some order-preserving map
$$t:\{0,1,...,2m+1\}\rightarrow\{0,1,...,2n+1\}$$
(the choice for such $t$ may not be unique).

First we can eliminate all zeros in $q_0,...,q_m$. Assume that $q_i=0$ for some $0\leq i\leq m$. Then the map
$$u:\{0,1,...,2m-1\}\rightarrow\{0,1,...,2m+1\}$$
defined by sending $0,1,...,2m-1$ to $0,1,...,2m+1$ in order, with $i$ and $2m+1-i$ removed, is symmetric and identifies $y$ with
$$y^\prime=(d_1,...,d_{i-1},d_id_{i+1},d_{i+2},...,d_m,d,\tau d_m,...,\tau d_{i+2},\tau(d_{i+1})\tau(d_i),\tau d_{i-1},...,\tau d_1)$$
$$\times (q_0,...,q_{i-1},q_{i+1},...,q_m,q_m,...,q_{i+1},q_{i-1},...,q_0)\in S_{m-1}$$
if $i<m$, or
$$y^\prime=(d_1,...,d_{m-1},d_md\tau d_m,\tau d_{m-1},...,\tau d_1)\times(q_0,...,q_{m-1},q_{m-1},...,q_0)\in S_{m-1}$$
if $i=m$. Since $u$ is symmetric, the images of $y$ and $y^\prime$ are identified. Since $x$ and $y^\prime$ are identified by the map $t\circ u$, it suffices to show that the images of $x$ and $y^\prime$ are identified. Hence we can replace $y$ by $y^\prime$. After finitely many such steps, we can eliminate all zeros inside $\{q_0,...,q_m\}$. So we can assume at the beginning that there are no zeros in $\{q_0,...,q_m\}$.

Assume that $t$ is not symmetric. Then there exists $0\leq i\leq m$ such that $t(j)+t(2m+1-j)=2n+1$ for $j=0,1,...,i-1$, but $t(i)+t(2m+1-i)\neq 2n+1$. Assume that $t(i)+t(2m+1-i)<2n+1$ (the proof for the other case is symmetric).

Since $t$ identifies these two points, we have
$$p_{t(i)}=\sum_{k\in t^{-1}(t(i))}q_k=\sum_{k<i,t(k)=t(i)}q_k+q_i+\sum_{k>i,t(k)=t(i)}q_k,$$
$$p_{2n+1-t(i)}=\sum_{l,t(l)=2n+1-t(i)}q_l.$$
They should be the same since $x$ is fixed. However, since $t(2m+1-i)<2n+1-t(i)$, $t(l)=2n+1-t(i)$ only if $l>2m+1-i$, in which case $t(2m+1-l)=2n+1-t(l)=t(i)$. Thus if $q_l$ appears in the sum
$$p_{2n+1-t(i)}=\sum_{l,t(l)=2n+1-t(i)}q_l,$$
$q_{2m+1-l}$ must appear in the first sum of
$$p_{t(i)}=\sum_{k<i,t(k)=t(i)}q_k+q_i+\sum_{k>i,t(k)=t(i)}q_k,$$
and we have $q_l=q_{2m+1-l}$ since $y$ is fixed. Thus
$$p_{t(i)}\geq\sum_{k<i,t(k)=t(i)}q_k+q_i\geq \sum_{l,t(l)=2n+1-t(i)}q_l+q_i=p_{2n+1-t(i)}+q_i>p_{2n+1-t(i)},$$
which is a contradiction. Thus $t$ must be symmetric, and hence the images of $x$ and $y$ are identified.

Therefore, the homeomorphisms between $S_n$'s and $T_n$'s pass through the quotient, and we get a homeomorphism $$\eta:(B^\sigma A)^{\mathcal{C}_2}\rightarrow B(*,A,A^{\mathcal{C}_2}).$$

\paragraph{(vii)} Now we can check that the map $H$ defined in Theorem \ref{mainlemma} gives a homotopy from $\iota$ to $p$.

The map
$$H:(B^\sigma A)^{\mathcal{C}_2}\times I\rightarrow B^\sigma A$$
given by sending
$$(a_1,...,a_n,a,\tau a_n,...,\tau a_1)\times (v_0,...,v_n,v_n,...,v_0)\times t$$
to
$$(a_1,...,a_n,a,\tau a_n,...,\tau a_1)\times((1+t)v_0,...,(1+t)v_n,(1-t)v_n,...,(1-t)v_0)$$
is well-defined since this definition is compatible with all symmetric relations, and we have shown in part (vi) that any pair of identified points can be related by a sequence of symmetric relations.

It's easy to check that $H(x,0)=\iota(x)$. And for any point 
$$(a_1,...,a_n,a,\tau a_n,...,\tau a_1)\times(v_0,v_1,...,v_n,v_n,...,v_0)\in(A^{2n+1}\times\Delta^{2n+1})^{\mathcal{C}_2},$$
its image through
$$(B^\sigma A)^{\mathcal{C}_2}\cong B(*,A,A^{\mathcal{C}_2})\xrightarrow{p} B(*,A,*)=B^\sigma A$$
is the point expressed by 
$$(a_1,...,a_n)\times(2v_0,2v_1,...,2v_n),$$
which is equivalent to the point 
$$(a_1,...,a_n,a,\tau a_n,...,\tau a_1)\times(2v_0,2v_1,...,2v_n,0,0,...,0)\in A^{2n+1}\times\Delta^{2n+1}.$$
Thus we also have $H(x,1)=p(\eta(x))$.

\bibliographystyle{abbrv}

\end{document}